\documentclass[11pt,a4paper]{article}
\usepackage{geometry}  
\usepackage[british]{babel} 
\usepackage{a4wide}
\usepackage[T1]{fontenc}
\usepackage[utf8]{inputenc}
\usepackage[pdftex]{graphicx}  
\usepackage{concrete,url} 
\usepackage{varioref}\labelformat{equation}{(#1)}
\usepackage{enumerate,mathrsfs,url}
\usepackage{amsfonts}
\newcommand{\anymode}[1]{\ifmmode{#1}\else\mbox{$#1$}\fi}
\newcommand{\biblioitem}[1]{\par \frenchspacing
           \vbox{\noindent \hangindent=2em \hangafter=1{#1}}
           \vskip\parskip}
\def\biblioitem{\par\vskip\parskip\noindent\hangindent=2em \hangafter=1}

\newcommand\NB[1]{\kern-0.4em\raisebox{-1.5ex}{$\stackrel{\big|}{\hbox{%
  \tiny\sc NB}}$}\marginpar{\sl\footnotesize#1\hfill}\kern-0.2em}

\long\def\Nota#1{\footnote{#1}\kern-0.2em\NB{cfr Nota N.\,\arabic{footnote}}}

\def\linea#1{\ifhmode\hfill\break\fi\hbox to \hsize{#1}}

\newcommand{\inv}{^{-1}}

\def\vc{v\kern -0.1em .c.\relax}

\newcommand{\dfrac}[2]{\displaystyle{\frac{#1}{#2}}}

\newcommand{\E}[2][]{
   \ensuremath{\mathbb{E}_{#1}\!\left\{\displaystyle{#2}\right\}}}

\newcommand{\half}{\mbox{$\textstyle \frac{1}{2}$}}   
\long\def\ignore#1{}

\newcommand{\indep}{\perp\kern-0.5em\perp}

\newcommand{\pr}[2][]{
   \ensuremath{\mathbb{P}_{#1}\!\left\{\displaystyle{#2}\right\}}}

\newcommand{\Real}{\mathbb{R}}

\newcommand{\T}{^{\top}}

\newtheorem{theorem}{Theorem}
\newtheorem{lemma}[theorem]{Lemma}

\newtheorem{corollary}[theorem]{Corollary}
\newtheorem{proposition}[theorem]{Proposition}

\newcommand{\SN}{\mathrm{SN}{}}

\newcommand{\N}{\mathrm{N}{}}

\renewcommand{\d}{\,\mathrm{d}}
\newcommand{\pd}[2]{\displaystyle\frac{\partial #1}{\partial #2}}

\title{Yet another skew-elliptical family but of a different kind:
   return to Lemma 1}
\author{
  {\Large Adelchi Azzalini} \\  \large 
  Dipartimento di Scienze Statistiche \\
  Università di Padova\\
  Italia 
  \and
  {\Large Giuliana Regoli} \\
   Dipartimento di Matematica e Informatica\\ 
    Università di Perugia\\
    Italia
 }
\date{\footnotesize\today}
\let\phi=\varphi

\begin{document}
\maketitle
\begin{abstract}
In the context of modulated-symmetry distributions, there exist various forms
of skew-elliptical families. We present yet another one, but with an unusual
feature: the modulation factor of the baseline elliptical density is
represented by a distribution function with an argument which is not an odd
function, as it occurs instead with the overwhelming majority of similar
formulations, not only with other skew-elliptical families. The proposal
is obtained by going back to the use of Lemma~1 of Azzalini and Capitanio
(1999), which can be seen as the general frame for a vast number of existing
formulations, and use it on a different route.
The broader target is to show that this `mother lemma' can still 
generate novel progeny.
\end{abstract}  

\vspace{2ex}

\noindent\emph{Some key-words:} elliptical distributions,
skew-elliptical distributions, symmetry-modulated distributions.
\clearpage

\section{Background and aims}

In the last fifteen years or so, there has been a formidable development 
on the theme of `symmetry-modulated distributions', also called
`skew-symmetric distributions'.  In the context of continuous random
variables, the key tool for building a \emph{symmetry-modulated distribution}
starting from a centrally symmetric density $f_0$ on $\Real^d$,
that is, such that $f_0(x)=f_0(-x)$,  is represented  by the expression
\begin{equation}
  f(x) = 2\: f_0(x) \:G_0\{w(x)\}\,,\qquad\quad x\in\Real^d,
  \label{e:SMD-pdf}
\end{equation}
where $G_0$ is a univariate continuous distribution function such that
$G_0(-x)=1-G_0(x)$, that is, it establishes to a symmetric distribution
about $0$, and $w(\cdot)$ is a real-valued function such that $w(-x)=-w(x)$.
Under the stated conditions, for any `baseline' density $f_0$ and 
any choice of the ingredients $G_0$ and $w$ which  enter the final
`modulating term' of $f(x)$, this is a proper density.

The above passage has reproduced Proposition~1 of Azzalini \& Capitanio (2003). 
A slightly different formulation, essentially equivalent to the one above, 
has been presented by Wang \emph{et alii} (2004).
Even a brief summary of the developments stemming from the symmetry-modulated 
distributions in \ref{e:SMD-pdf} would take an enormous space. 
For such objective, the
reader is referred to the recent account by Azzalini \& Capitanio (2014). 
Here we limit ourselves to recall the facts directly related to our plan of work.

The fact that \ref{e:SMD-pdf} constitutes a density function is actually
a corollary of a more general result, namely Lemma~1 of Azzalini \& 
Capitanio (1999); its exact statement will be recalled in the next section.
Right after this lemma, Azzalini and Capitanio formulated an immediate 
corollary by taking $f_0$ of elliptical type and $w$ of linear type.
That corollary prompted a stream of literature focusing on successive
layers of generalizations, where $f_0$ fulfils some symmetry condition 
and $w$ is odd, eventually leading to symmetry-modulated distributions.

There exist, however, other constructions encompassed by the above-quoted
Lemma~1, falling outside the domain of symmetry-modulated distributions
in \ref{e:SMD-pdf}.  This territory has been barely explored. 
As far as we known, the only two publications of this sort are 
those of  Azzalini (2012) and Jupp \emph{et alii} (2016).
The broad target of the present note is to investigate further this area.
In spite of the common domain of interest, operationally there is 
little overlap between our development and the publications just quoted.

The present specific contribution, developed along the indicated direction, 
is constituted by a new proposal of skew-elliptical  or, 
equivalently, skew-elliptically contoured (SEC) density.
A number of SEC constructions in the form of  symmetry modulated densities, 
obtained by some form of perturbation  of an elliptically contoured (EC) 
density, have already been examined in the literature.
Besides the above-mentioned `linear-type' SEC, the more commonly 
adopted form of SEC is the one introduced by Branco \& Dey (2001)
via a selection mechanism applied to an elliptical distribution. 
In their exposition, it was not evident that their SEC distribution had
a structure like \ref{e:SMD-pdf}, but this was shown in some important
special cases by Azzalini \& Capitanio (2003) and later in general
by Azzalini \& Regoli (2012, p.\,872).
An additional form of SEC distribution has been presented by Sahu 
\emph{et alii} (2003).

There is structural difference between the existing SEC distributions 
and  the one to be presented here.
Because of its genesis, the new SEC is not constrained by the conditions
pertaining to \ref{e:SMD-pdf}. It obviously retains the fact
that $f_0$ is of elliptical class, otherwise the term SEC would not
apply to it, but its $w(\cdot)$ does not satisfy the condition 
$w(-x)=-w(x)$.


\section{A `mother' lemma}

The wording of the next statement is slightly different from the one 
in the original source, but completely equivalent to it.
The expression `symmetric about $0$' is an abbreviated term for
`symmetrically distributed about $0$', in case of a univariate
random variable. We also use the term `symmetric about $0$' for
a univariate distribution function $H$ such that $H(-x)=1-H(x)$.

\begin{lemma}[Azzalini \& Capitanio, 1999] \label{th:mother-lemma} 
  Denote by $G_0$  the continuous distribution function of a univariate 
  random variable symmetric about $0$, by $w(\cdot)$ 
  a real-valued function on $\Real^d$ and by $Z_0$ a $d$-dimensional 
  variable with density function $f_0$, such that  $W=w(Z_0)$ 
  is symmetric about $0$. Then
  \begin{equation} \label{e:mother-lemma1}
        f(x) = 2\:f_0(x)\:G_0\{w(x)\} \,,     \qquad x\in\Real^d,
  \end{equation}  
  is a density function.
\end{lemma}

At first sight, the distribution in \ref {e:mother-lemma1} coincides with
\ref{e:SMD-pdf}. While the mathematical expressions look the same, the
meaning of their symbols, and consequently their mathematical meaning, 
are largely different:
\begin{itemize}
\item in \ref {e:mother-lemma1}, $f_0$ is not required to satisfy the 
   condition $f_0(-x)=f_0(x)$;
\item    in \ref {e:mother-lemma1}, $w$ is not required to satisfy the 
   condition $w(-x)=-w(x)$;
\item the only common assumption is that $G_0(-x)=1-G_0(x)$. 
\end{itemize}

The key point for applying Lemma~\ref{th:mother-lemma} is that $W=w(Z_0)$
is symmetric about 0.
Once this holds, any  $G_0$ symmetric about 0 can be adopted to build a valid 
density. This scheme was particularly simple in the first application
of the lemma, namely Corollary~2 of Azzalini \& Capitanio (1999):
if $Z_0$ is elliptical centred at 0, then the linear transform $w(Z_0)=a\T\,Z_0$ 
is symmetric about 0 for any choice of the vector $a$ of coefficients 
and this suffices for applying the lemma,
irrespective of the choice of $G_0$. It only takes a little more effort
along the same line of argument to show that $f(x)$ in \ref{e:SMD-pdf} 
is a proper density.

We also recall the stochastic representation for distribution 
\ref{e:mother-lemma1} stated by Azzalini \& Capitanio (1999, p.\,599). 
If $X$ is a $d$-dimensional  variable with density $f_0$ 
and $T$ is an independent variable with distribution function $G_0$, 
satisfying the conditions of Lemma~\ref{th:mother-lemma}, then
\begin{equation}
  Z = \cases{X   & if $T \le w(X)$ \cr
             -X  & otherwise}
  \label{e:mother-lemma-SR}                        
\end{equation}
has distribution \ref{e:mother-lemma1}.

Besides symmetry-modulated distributions, Lemma~\ref{th:mother-lemma}
embraces other constructions, although so far these have little been 
studied, as already mentioned. 
This overarching capability of Lemma~\ref{th:mother-lemma} explains
the term `mother lemma' which we have employed.
The aim of the present note is to examine a case of this sort, specifically
when $f_0$ is an elliptical density.


\section{Functions of orthogonal elliptical components}

We want to examine a construction where the density function $f_0$ of $Z_0$ is
of EC type and $W=w(Z_0)$ is symmetric about 0.
 
Start by considering the case where $f_0$ is a bivariate elliptically contoured
distribution with `standardized' marginals. 
Specifically, assume that $f_0$ is the density function of
\begin{equation}
  \pmatrix{X \cr Y} \sim 
  \mathrm{EC}_2\left(\matrix{0 \cr 0}, \pmatrix{1 & \rho \cr \rho & 1}, \psi \right)
  \label{e:std-EC2}
\end{equation}
following the notation adopted by Fang \emph{et alii} (1990), up to a minor change 
in the symbol of the `generator', $\psi$ for us.

\begin{proposition} 
\label{th:symm-W}
Given the continuous random variable $(X,Y)$ with distribution as 
in \ref{e:std-EC2}, the transformed variable 
\begin{equation}
   W = (Y- \rho X) \: h(X)    \label{e:W-simple}
\end{equation}
has is  symmetric about 0 for any choice of the  real-valued function $h$.
\end{proposition}
\emph{Proof}
From Theorem 2.18 of Fang \emph{et alii} (2012), the conditional distribution of $W$ given $X=x$ is
\[
   (W|X=x) \sim \mathrm{EC}_1(0, (1-\rho^2)\,h(x)^2, \psi_x)
\]
where $\psi_x$ is another generator, which in general depends on $x$, 
but its explicit expression is irrelevant to us.
All we need is that the corresponding density $f_{W|x}(w)$ is symmetric about~0, 
so that $f_{W|x}(w)=f_{W|x}(-w)$. Therefore the unconditional density of $W$ is
\begin{eqnarray*}
  f_W(w) &=& \int_\Real f_{W|x}(w)\: f_{0X}(x) \d{x} \\
         &=& \int_\Real f_{W|x}(-w)\: f_{0X}(x) \d{x} \\
         &=& f_W(-w)\,. 
\end{eqnarray*}
by integration with respect to the distribution $f_{0X}(x)$ of the $X$ component.
\hfill \textsc{qed}

This result can  immediately be extended to a general $d$-dimensional continuous
random variable  with EC distribution $(d>1)$, whose components are partitioned 
as follows:
\begin{equation}
  \pmatrix{X \cr Y} \sim 
  \mathrm{EC}_d\left(\mu, \Sigma, \psi \right),
  \qquad \mu= \pmatrix{\mu_x \cr \mu_y},
  \quad \Sigma =\pmatrix{\Sigma_{11} & \sigma_{12} \cr \sigma_{21} & \sigma_{22}}
  \label{e:EC}
\end{equation}
where now $X$ is $(d-1)$-dimensional; $W$ takes now the more general form
\begin{equation}
   W = \big(Y- m_Y(X)\big) \: h(X)  \, \quad \mathrm{or} \quad
   W = \frac{Y- m_Y(X)}{(\sigma_{22\cdot1})^{1/2}} \: h(X)  
   \label{e:W}
\end{equation}
where
\begin{eqnarray}
 m_Y(X) &=& \mu_y + \sigma_{21}\Sigma_{11}\inv(X-\mu_x) = \beta_0 + \beta\T X 
     \,,~ \mathrm{say},  \label{e:E(Y|x)}\\  
 \sigma_{22\cdot1} &=& \sigma_{22} -\sigma_{21}\Sigma_{11}\inv\sigma_{12}
     \label{e:sigma22.1}
\end{eqnarray}
are as in equation (2.42) of Fang \emph{et alii} (2012). 
In case the mean value of $Y$ conditionally on $X$, 
it coincides with  \ref{e:E(Y|x)}.
The second form of \ref{e:W} is just a minor variant of the first one,
but it is conceptually convenient to have the first factor `standardized',
in the sense of being free from scale factors, besides the location parameter.
Finally, note that symmetry about 0 persists if $Y- m_Y(X)$ is transformed
to $w_0(Y- m_Y(X))$ by an odd function $w_0$.

We can summarize the combination of the above discussion and  
Lemma~\ref{th:mother-lemma} into the following conclusion.

\begin{proposition} \label{th:main}
Denote by $f_0(x,y)$ the $d$-dimensional elliptically contoured density 
of the random variable $(X,Y)$ in \ref{e:EC}, by $G_0$ a continuous distribution
function such that $G_0(-x)=1-G_0(x)$ and by $w(x,y)$ either of 
the two forms
\begin{equation}
   w(x,y) = w_0\big(y- m_Y(x)\big) \: h(x)  \,, \quad \quad
   w(x,y) = w_0\left(\frac{y- m_Y(x)}{(\sigma_{22\cdot1})^{1/2}}\right)\: h(x)
   \label{e:w}
\end{equation}
for any function $h(x)$ from $\Real^{d-1}$ to $\Real$ and any odd function 
$w_0$ on the real line; here $m_Y(\cdot)$ and $\sigma_{22\cdot1}$
are given by \ref{e:E(Y|x)} and \ref{e:sigma22.1}.
Then 
\begin{equation}
    f(x,y) = 2\:f_0(x,y) \:G_0\{w(x,y)\}
    \label{e:f0-modulated}
\end{equation}
is a proper density function.
\end{proposition}

\begin{figure}
\centering
{
\includegraphics[width=0.45\hsize]{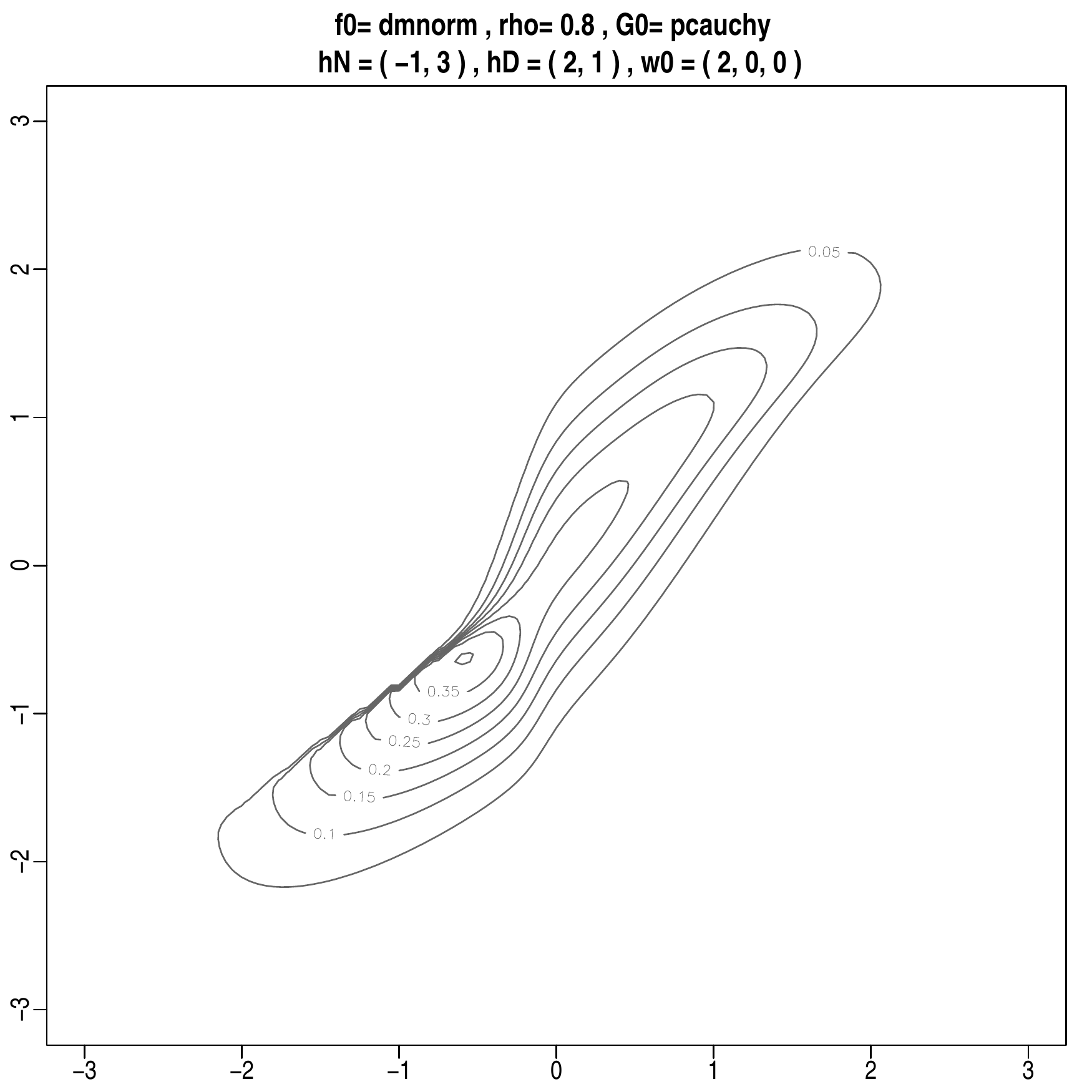}
\includegraphics[width=0.45\hsize]{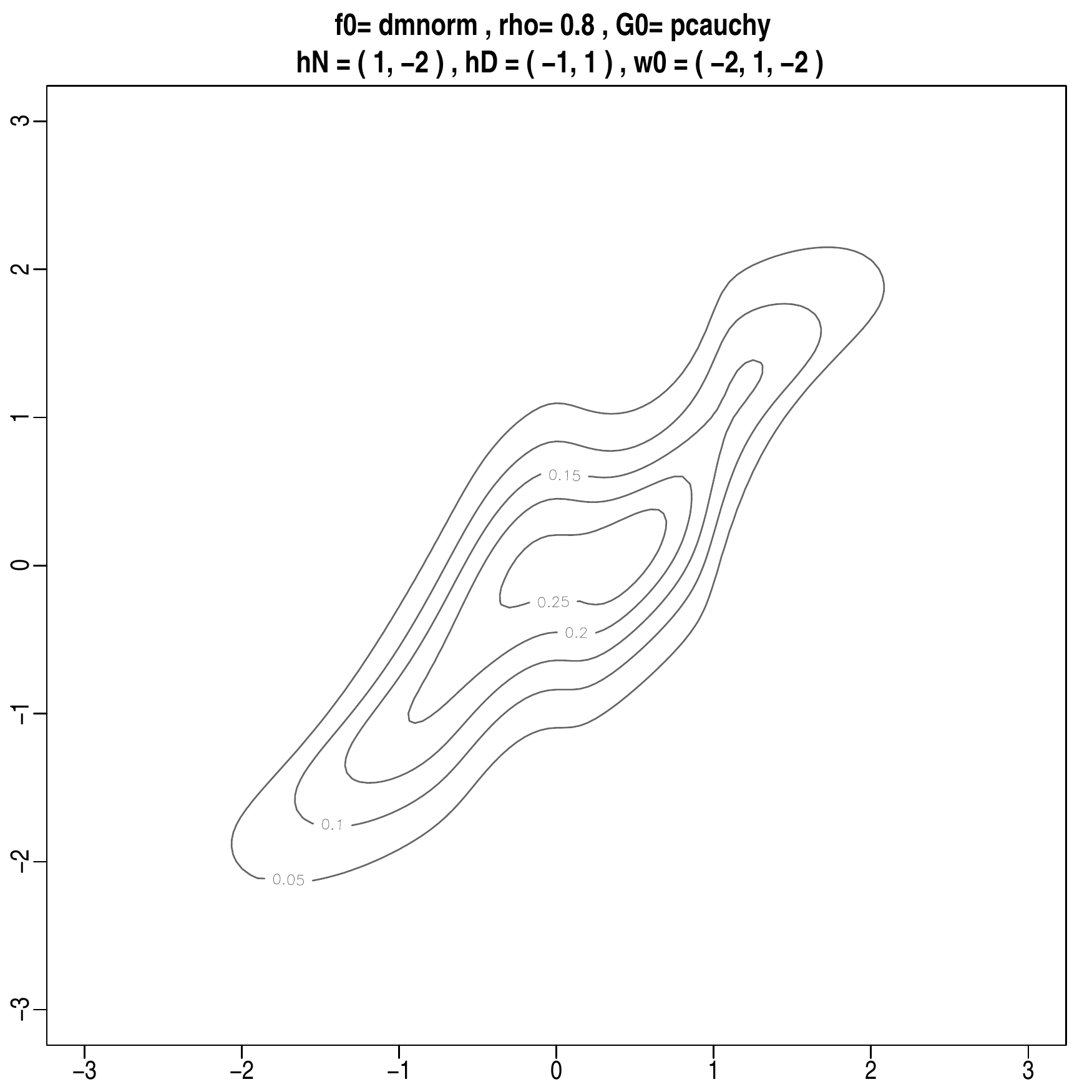}
}
\par
{
\includegraphics[width=0.45\hsize]{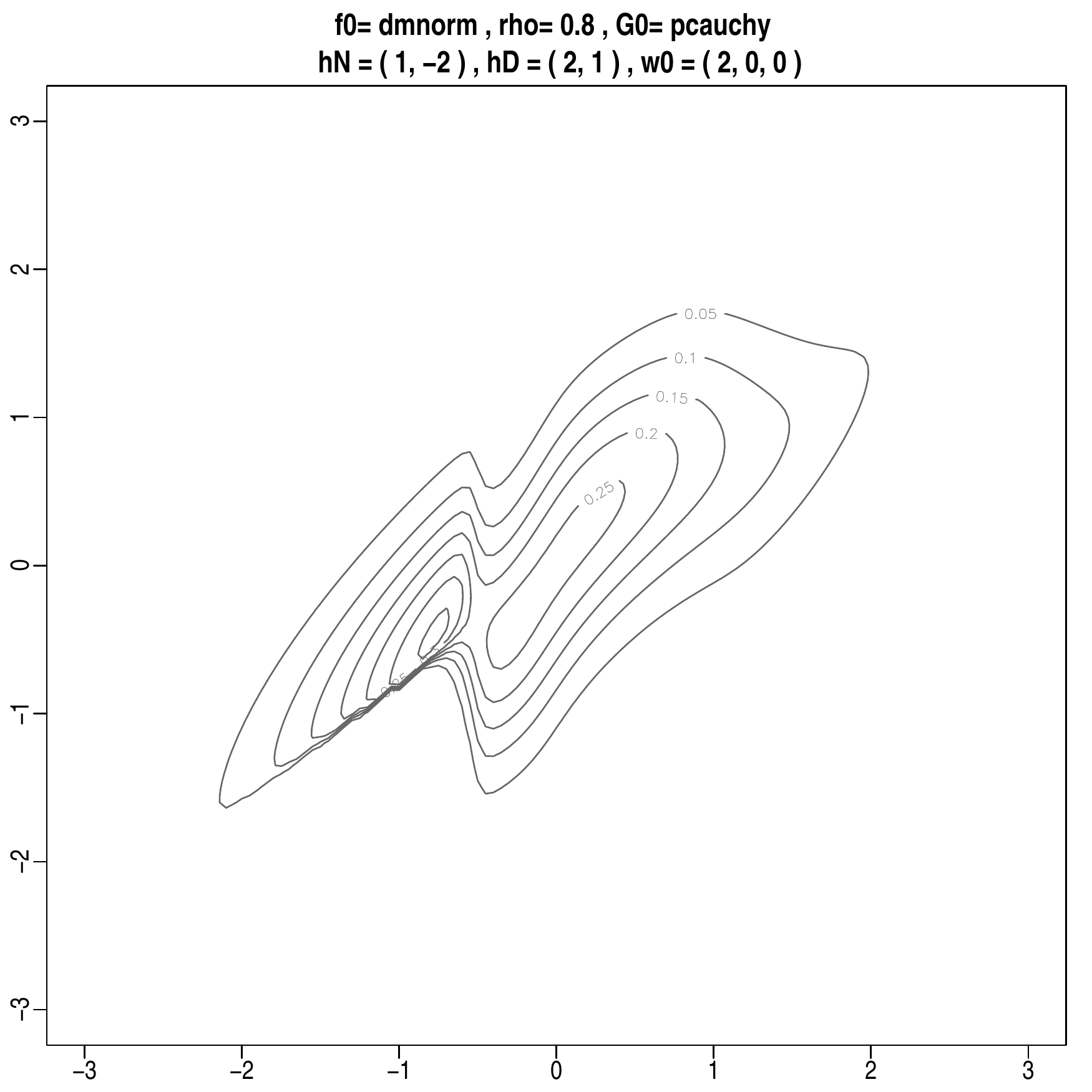}
\includegraphics[width=0.45\hsize]{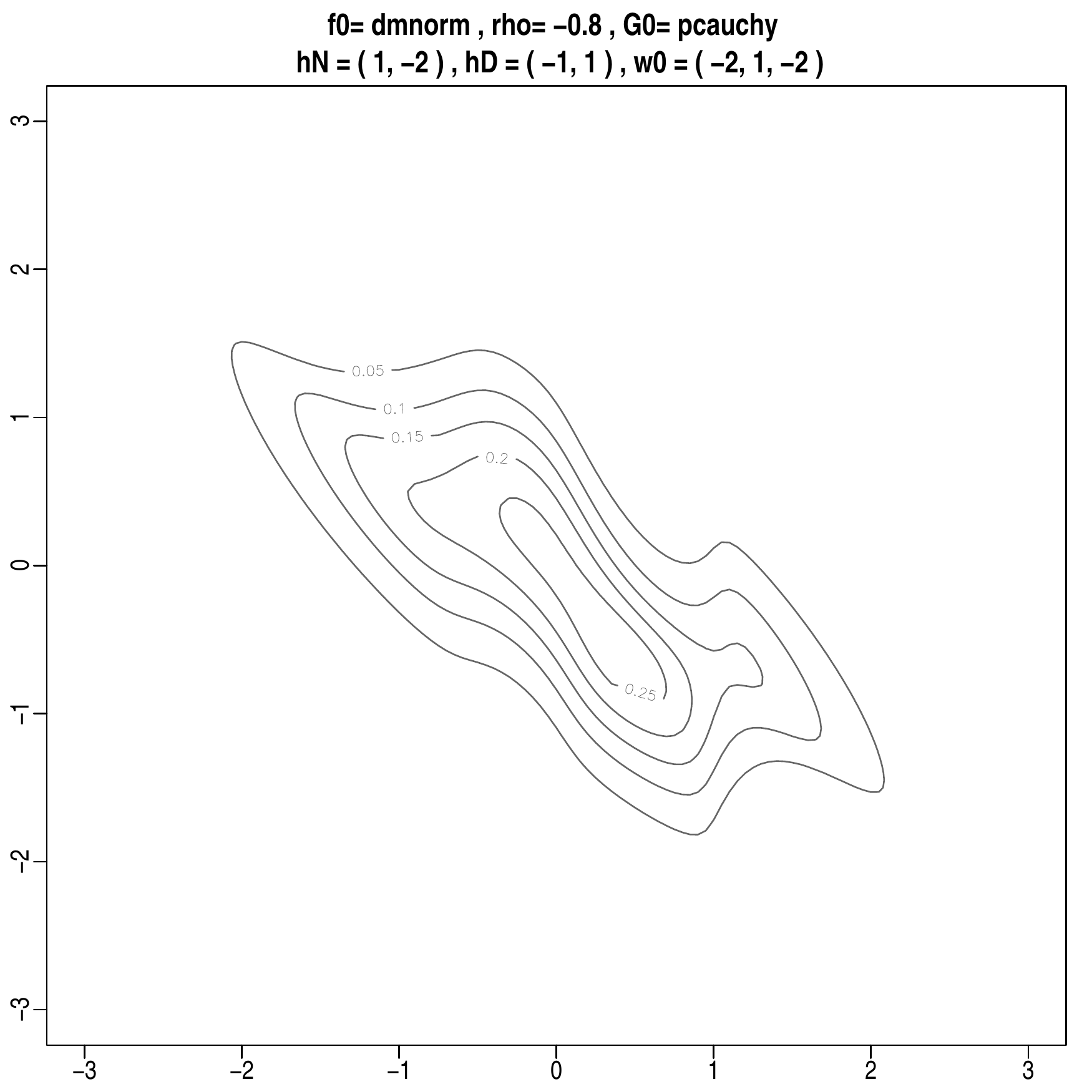}

}
\par{
\includegraphics[width=0.45\hsize]{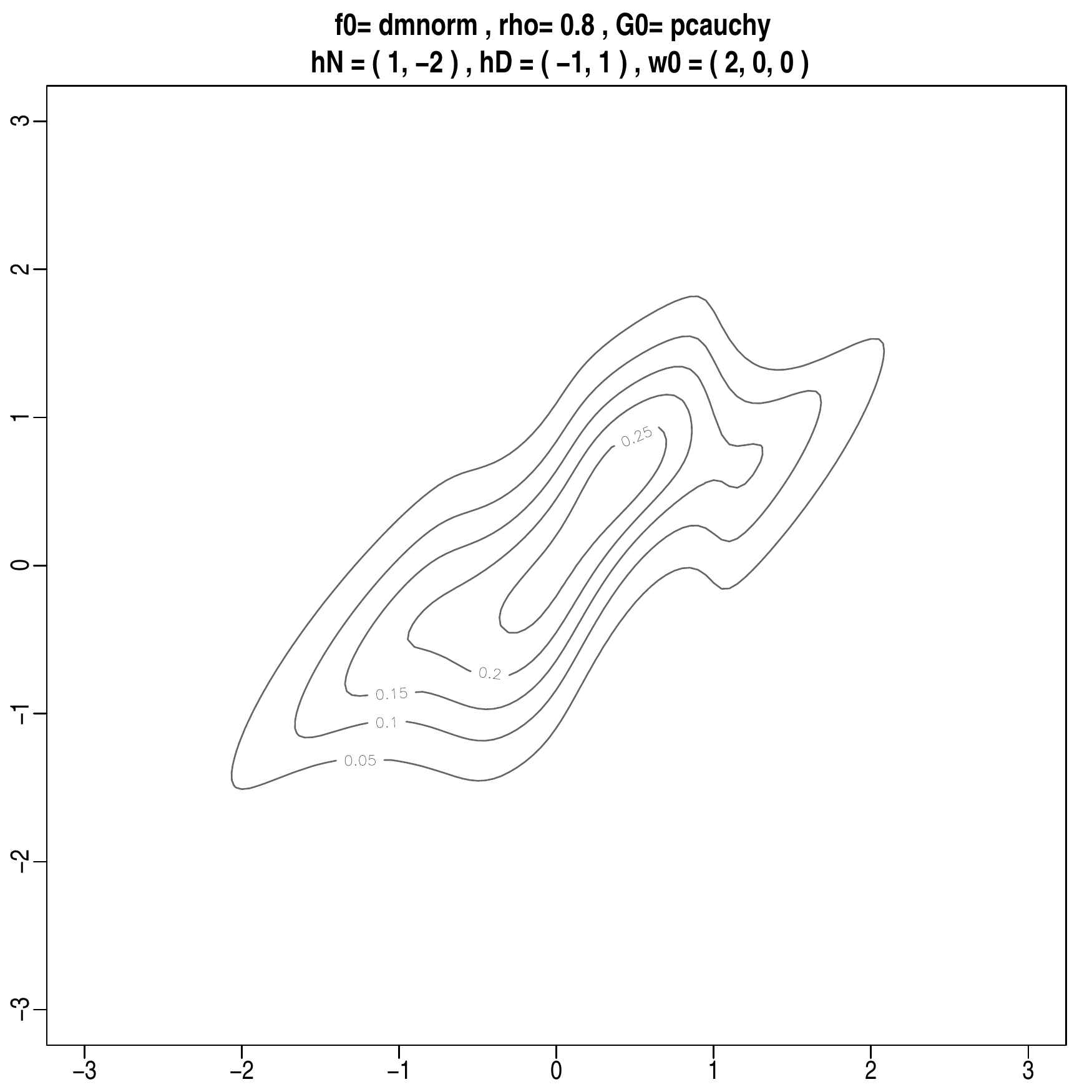}
\includegraphics[width=0.45\hsize]{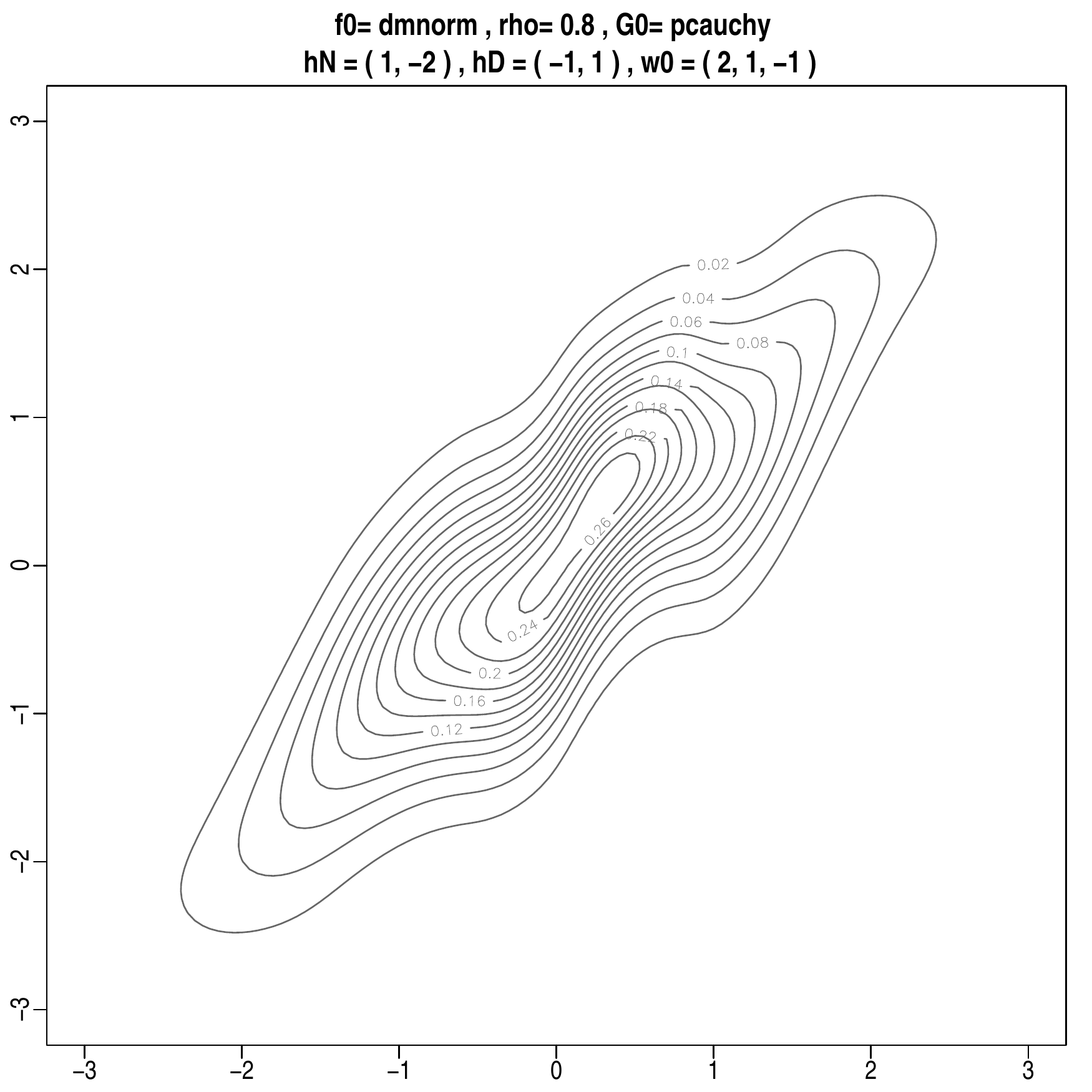}

}
\caption{Example plots of densities defined in Proposition~\ref{th:main}
using functions in \ref{e:hN,hD,w0}.}
\label{f:density-examples}  
\end{figure}

Some example bivariate densities $f(x,y)$ of the type introduced in 
Proposition~\ref{th:main} are displayed in Figure~\ref{f:density-examples} 
in the form of contour level plots. The second form of \ref{e:w} has been
used and the other ingredients are as follows: 
 $f_0(x)=\phi_2(x; \Sigma)$  where $\Sigma$ is the same matrix
occurring in \ref{e:std-EC2}, $G_0$ is the standard Cauchy distribution 
function and
\begin{equation}
  h(x)= \frac{h_N(x)}{h_D(x)} =\frac{1+a_1\,x+a_2\,x^2}{1+b_1\,x+b_2\,x^2}\,,
  \hspace{4em}  
  w_0(x) = \frac{c_1\,x + c_3\,x^3}{1+ c_2\,x^2}
  \label{e:hN,hD,w0}
\end{equation}
for various choices of the coefficients $a, b, c$ and the correlation $\rho$ 
in $\Sigma$. The specific choices of the coefficients  $a, b, c$ are 
indicated at the top of each pane of Figure~\ref{f:density-examples}.
The plots indicate a wide flexibility of the family of distribution, even
employing a relatively limited number of cofficients in the ratios
\ref{e:hN,hD,w0}.


The general stochastic representation \ref{e:mother-lemma-SR} applies directly 
to distribution \ref{e:f0-modulated}.

\section{The unperturbed component} \label{s:X-marginal}

Consider the $d$-dimensional density
\begin{equation}
   f(x,y) = 2\: f_0(x,y)\: G_0\{(y-\beta_0 -\beta\T x)\: h(x)\},
   \qquad\qquad x\in\Real^{d-1},\quad y\in\Real,
   \label{e:f(x,y)}
\end{equation}
where $f_0$ is as in \ref{e:EC} and we use the expression in \ref{e:E(Y|x)}.
If $(\tilde{X}, \tilde{Y})$ is a random variable with density \ref{e:f(x,y)},
the marginal distribution of $\tilde{X}$ is
\begin{eqnarray*}
  f_{\tilde X}(x) &=& \int_{\Real} f(x,y)\d{y} \\
         &=& \int_{\Real} 2\: f_{0,X}(x)\: f_{0,Y|x}(y)\:
              G_0\{(y-\beta_0 -\beta\T x)\: h(x)\}\d{y}\\
         &=&  f_{0,X}(x) \int_{\Real} 
              2\:f_{0,Y|x}(z)\: G_0(z\:\sigma_{22\cdot1}\: h(x)) \d{z}\\
          &=&  f_{0,X}(x)     
\end{eqnarray*}
on recalling (i)~the expressions of the conditional mean \ref{e:E(Y|x)} and
the scale factor \ref{e:sigma22.1}, (ii)~the fact that $f_{0,Y|x}(\cdot)$ is 
symmetric about its mean (conditional) value, (iii)~Lemma~1 of Azzalini (1985).

Therefore the marginal distribution $f_{\tilde X}$ of the first ${d{-}1}$ 
component variables after the perturbation operation is equal to the original 
unperturbed density $f_{0,X}$.


\section{Moment generating function in a simple case}

We want to compute the moment generating function (MGF) for the basic case 
\ref{e:std-EC2} of normal type, $W$ as in \ref{e:W-simple} and $G_0=\Phi$. 
In other words, consider the density
\begin{equation}
   f(x,y) = 2 \phi(x, y; \rho) \Phi\{(y-\rho x)\:h(x)\}\,, 
   \qquad (x,y)\in\Real^2.
   \label{e:norm2perturb}
\end{equation}
For notational convenience, we modify slightly the notation and, from 
now on, we use $(X,Y)$ to denote a random variable with distribution 
\ref{e:norm2perturb}.
The MGF of $(X,Y)$  is
\begin{eqnarray*}
  M(t_1, t_2) 
    &=& \int_\Real \int_\Real \exp(t_1 x+ t_2 y) \:f(x,y) \d{x}\d{y}\\
    &=& 2 \int_\Real e^{t_1 x} M_c(t_2, x) \phi(x) \d{x}
\end{eqnarray*}
where
\begin{eqnarray*}
  M_c(t_2, x) 
    &=& \int_\Real e^{t_2 y}\: f_{0,Y|x}(y) \:\Phi\{(y-\rho x)\:h(x)\} \d{y} \\
    &=& \int_\Real e^{t_2 y} \phi\left(\frac{y-\rho x}{\sqrt{1-\rho^2}}\right)
          \frac{1}{\sqrt{1-\rho^2}} \:\Phi\{(y-\rho x)\:h(x)\} \d{y} \\
    &=& \exp(t_2\rho x)  \int_\Real \exp\left(t_2 z \sqrt{1-\rho^2}\right) \:
           \phi(z)\:\:\Phi\left(z\sqrt{1-\rho^2}\:h(x)\right) \d{z} \\
    &=& \exp(t_2\rho x) \int_\Real \frac{1}{\sqrt{2\pi}}\:
        \exp\left\{-\half\left[z^2 -2\, t_2 z \sqrt{1-\rho^2}
        \pm t_2^2(1-\rho^2)\right]\right\} \:
             \Phi\left(z\sqrt{1-\rho^2}\:h(x)\right) \d{z} \\
    &=& \exp\left[t_2\rho x +\half t_2^2(1-\rho^2)\right] \int_\Real \phi(u)\,
        \Phi\left\{\left(u+ t_2\sqrt{1-\rho^2}\right) \:h(x)\right\} \d{u}  \\
    &=& \exp\left[t_2\rho x +\half t_2^2(1-\rho^2)\right]  \:
        \Phi\left(\frac{t_2\, (1-\rho^2)\,h(x)}{\sqrt{1+(1-\rho^2)\,h(x)^2}}\right)         
\end{eqnarray*}
where the last equality follows from Corollary 1 of Ellison's (1964) Theorem 2.
Therefore
\begin{eqnarray*}
  M(t_1, t_2) 
    = 2\:\exp\left\{\half t_2^2(1-\rho^2)\right\}\: 
     \int_\Real \exp\{x(t_1+t_2\rho)\}\: 
      \Phi\left(\frac{t_2\, (1-\rho^2)\,h(x)}{\sqrt{1+(1-\rho^2)\,h(x)^2}}\right)    
      \phi(x) \d{x}\,.
\end{eqnarray*}     
The above integral does not lend itself to explicit solution. To compute
the moment of the $X$ and $Y$ components of \ref{e:norm2perturb}, we exchange
the integration and differentiation steps; hence consider
\[
   \left.\pd{M}{t_1}\right|_{t_1=t_2=0} \quad \mathrm{and}\quad 
   \left.\pd{M}{t_2}\right|_{t_1=t_2=0}
\] 
with differentiation under the integration sign.
The first expression leads to $\E{X}=0$, as expected, since the marginal density of
$X$ is $\phi(x)$, from \S\,\ref{s:X-marginal}. The second expression leads to
\begin{equation}
  \E{Y}= \sqrt{\frac{2}{\pi}}\:(1-\rho^2)\: 
         \int_\Real \frac{h(x)}{\sqrt{1+(1-\rho^2)\,h(x)^2}}\:\phi(x)\d{x} \,.
   \label{e:norm2-EY}      
\end{equation}

\paragraph{A very special case} As a simple check of the above result, 
consider the elementary case $h(x)\equiv 1$,
for which \ref{e:norm2-EY} lends 
\[    \sqrt{\frac{2}{\pi}}\:\frac{(1-\rho^2)}{\sqrt{2-\rho^2}}\,. \] 
In this case,  \ref{e:norm2perturb} is of type 
\[
  \SN_2\left(\xi=\pmatrix{0 \cr 0}, \Omega=\pmatrix{1 & \rho \cr \rho & 1},
               \alpha=\pmatrix{-\rho \cr 1} \right)
\]
and known formulae such as (5.31) and (5.11) of the Azzalini \& Capitanio (2014) 
say that
\[
   \E{\pmatrix{X\cr Y}}= \sqrt{\frac{2}{\pi}} \:
                         \pmatrix{0 \cr \dfrac{1-\rho^2}{\sqrt{2-\rho^2}}}.
\]
Addition numerical checks with other choices of $h$ confirm expression 
\ref{e:norm2-EY} by use of numerical integration.

\paragraph{Other feasible solutions}
The explicit expression of the integral in \ref{e:norm2-EY} is feasible 
only in favourable cases. 

A relatively common situation occurs when $h(x)$ is odd, 
so that the integrand  of \ref{e:norm2-EY} is odd and 
\[  \E{Y}=0 \,.\]
 Note that, when $h$ is odd,
$w(x,y)=(y-\rho x)\:h(x)$ is  an even function of $(x,y)$.
The special case $h(x)=\alpha x$ reduces to the distribution in
the first expression of equation (17) of Azzalini (2012).

Another class of functions for which explicit integration is feasible is
as follows. Suppose $s(x)$ is a function for which $\E{s(Z)}$ is known to be
$S$, say, when $Z\sim \N(0,1)$. Then solve 
\[
   \frac{h(x)}{\sqrt{1+ (1-\rho^2)\:h(x)^2}} = s(x)
\]
and obtain
\[
   h(x) = \frac{s(x)}{\sqrt{1 - (1-\rho^2)\:s(x)^2}} \,.
\]
For this choice of $h(x)$, we can say that,  by construction,
\[  \E{Y}= \sqrt{\frac{2}{\pi}}\:(1-\rho^2)\: S\,.\]
An especially simple example is $s(x)=\cos x$, such that $S=e^{-1/2}$,
leading to
\[
   h(x) = \frac{\cos x}{\sqrt{1-(1-\rho^2) \:(\cos x)^2}}\,,\qquad
   \E{Y}= \sqrt{\frac{2}{\pi\: e}}\:(1-\rho^2)\, .
\]

Another special case leading to a solution, with help from the Maxima
symbolic manipulation system, is $h(x)=\alpha |x|$ for which we obtain 
\[
  \E{Y} = \frac{2}{\sqrt{1-\rho^2}}\:
         \left[1-\Phi\left(\frac{1}{\alpha\sqrt{1-\rho^2}}\right)\right]
         \exp\left(\frac{1}{2\:\alpha^2 (1-\rho^2)}\right)\,.
\]

Although in many other cases \ref{e:norm2-EY} does not lend itself to
an explicit expression, it still simplifies computation with respect to 
direct two-dimensional  numerical  integration of $y\,f(x,y)$, since 
\ref{e:norm2-EY} requires only one-dimensional integration.

\section{Extension to multivariate $Y$}  \label{s:multiv-Y}

Consider the more general case where $Y$ is $m$-dimensional, with $1\le m<d$.
Hence write the joint distribution of $(X,Y)$ as
\begin{equation}
  \pmatrix{X \cr Y} \sim 
  \mathrm{EC}_d\left(\mu, \Sigma, \psi \right),
  \qquad \mu= \pmatrix{\mu_x \cr \mu_y},
  \quad \Sigma =\pmatrix{\Sigma_{11} & \Sigma_{12} \cr \Sigma_{21} & \Sigma_{22}}
  \label{e:EC_d}
\end{equation}
where $X$ is $(d-m)$-dimensional. 
We want to extend the expressions in \ref{e:W-simple} and \ref{e:W} to 
this new case.

\begin{proposition}\label{th:Y-multidim}
For a $d$-dimensional variable $(X,Y)$ distributed as in \ref{e:EC_d},
consider the transform 
\begin{equation}
  W = w_0\left(Y-m_Y(X)\right) \: h(X) 
  \label{e:w0*h}
\end{equation}
where $h$ is an arbitrary function function and $w_0$ is an odd function, 
that is, $w_0(u)\in\Real$ and $w_0(-u)=-w_0(u)$ for $u\in\Real^m$.
Then $W$ is distributed symmetrically around 0.
\end{proposition}
\emph{Proof.} Write $U=Y-m_Y(X)$ and $W_0=w_0(U)$. First we show that 
the distribution of $W_0$ conditionally on $X=x$ is symmetric about 0,
for any given $x$; in fact
\begin{eqnarray*}
  \pr{W_0 \le t |X=x} &=& \pr{-W_0 \ge -t |X=x} \\
     &=& \pr{w_0(-U) \ge -t |X=x} \\
     &=& \pr{w_0(U) \ge -t |X=x}  \\ 
     &=& \pr{W_0 \ge -t |X=x} 
\end{eqnarray*}
for any real $t$, and the last equality implies symmetry of $(W_0|X=x)$. 
Then argue like in Proposition~\ref{th:symm-W} to conclude that $W$ is
symmetric about 0.\hfill \textsc{qed} 

\begin{corollary}\label{th:coroll3}
Under the assumption of Proposition~\ref{th:main}, denote by $f_0$ 
the density of \ref{e:EC_d} and let $\beta_0+\beta\T x= m_Y(x)$ 
as given in equation (2.42) of Fang \emph{et alii} (2012).
Then expression \ref{e:f0-modulated} with 
\[ 
   w(x,y) = w_0( y- \beta_0 - \beta\T x)\: h(x)
\]
is a proper density function.
\end{corollary}

\section*{References}

 
\biblioitem
  Azzalini, A. (1985).
  {}A class of distributions which includes the normal ones.
{}{\em Scand.\ J.\ Statist.}, {\bf 12}, 171--178.

\biblioitem
  Azzalini, A. (2012).
  Selection models under generalized symmetry settings.
  \emph{Annals of the Institute of Statistical Mathematics},
   64, 737--750.

\biblioitem
   Azzalini, A., and Capitanio, A. (1999).
{}Statistical applications of the multivariate skew normal
  distribution.
{}{\em J.\ R.\ Statist.\ Soc., ser.\,B}, {\bf 61}, 579--602.
{}Full version of the paper at \url{arXiv.org:0911.2093}.
 
\biblioitem
 Azzalini, A., and Capitanio, A. (2003).
{}Distributions generated by perturbation of symmetry with emphasis on
  a multivariate skew $t$ distribution.
{}{\em J.\ R.\ Statist.\ Soc., ser.\,B}, {\bf 65}, 367--389.
{}Full version of the paper at \url{arXiv.org:0911.2342}.

\biblioitem  
A. Azzalini, with the collaboration of A. Capitanio (2014).
{\em The Skew-Normal and Related Families}. 
Cambridge University Press, Cambridge.

\biblioitem  
Azzalini, A., and Regoli, G. (2012).
{}Some properties of skew-symmetric distributions.
{}{\em Ann.\ Inst.\ Statist.\ Math.}, {\bf 64}, 857--879.
{}Available online 09 Sept 2011.

\biblioitem  
Ellison, B.~E. (1964).
{}Two theorems for inferences about the normal distribution with
  applications in acceptance sampling.
{}{\em J.\ Amer.\ Statist.\ Assoc.}, {\bf 59}, 89--95.

\biblioitem  
Fang, K.-T., Kotz, S., and Ng, K.~W. (1990).
{}{\em Symmetric multivariate and related distributions}.
{}London: Chapman \& Hall.

\biblioitem 
Jupp, P.\,E., Regoli, G. and Azzalini, A. (2016)
A general setting for symmetric distributions and their relationship 
to general distributions.
{\em Journal of Multivariate Analysis} 148, 107-–119.


\biblioitem
Sahu, K., Dey, D.~K., and Branco, M.~D. 2003.
{}A new class of multivariate skew distributions with applications to
  {B}ayesian regression models.
{}{\em Canad. J. Statist.}, {\bf 31}, 129--150.
{}Corrigendum: vol.\,37 (2009) , 301--302.

\biblioitem
 Wang, J., Boyer, J., and Genton, M.~G. (2004).
{}A skew-symmetric representation of multivariate distributions.
{}{\em Statist.\ Sinica}, {\bf 14}, 1259--1270.

\end{document}